\documentclass{amsart}
\usepackage{amsmath}
\usepackage{graphicx}
\usepackage{latexsym}

\newtheorem{thm}{Theorem}
\newtheorem{lem}[thm]{Lemma}
\newtheorem{cor}[thm]{Corollary}

\theoremstyle{definition}

\newtheorem{rmk}[thm]{Remark}

\newcommand{\CPb}{\overline{\mathbb{CP}}{}^{2}}
\newcommand{\CP}{{\mathbb{CP}}{}^{2}}
\newcommand{\R}{\mathbb{R}}
\newcommand{\Z}{\mathbb{Z}}
\newcommand{\SW}{{\rm SW}}

\def \x {\times}

\begin{document}

\title[Constructing infinitely many smooth structures]
{Constructing infinitely many smooth\\
structures on small 4-manifolds}

\author[A. Akhmedov]{Anar Akhmedov}
\address{School of Mathematics, Georgia Institute of Technology \newline
\hspace*{.375in} Atlanta, GA, 30332-0160, USA}
\email{ahmadov@math.gatech.edu}

\author[R. \.{I}. Baykur]{R. \.{I}nan\c{c} Baykur}
\address{Department of Mathematics, Michigan State University \newline
\hspace*{.375in} East Lansing, MI, 48824, USA}
\email{baykur@msu.edu}

\author[B. D. Park]{B. Doug Park}
\address{Department of Pure Mathematics, University of Waterloo \newline
\hspace*{.375in} Waterloo, ON, N2L 3G1, Canada}
\email{bdpark@math.uwaterloo.ca}

\date{February 2, 2007. Revised on April 15, 2007}


\begin{abstract}
The purpose of this article is twofold. First we outline a general construction scheme for producing simply-connected minimal symplectic $4$-manifolds with small Euler characteristics. Using this scheme, we illustrate how to obtain irreducible symplectic $4$-manifolds homeomorphic but not diffeomorphic to $\CP\#(2k+1)\CPb$ for $k = 1,\dots,4$, or to $3\CP\# (2l+3)\CPb$ for $l =1,\dots,6$.  Secondly, for each of these homeomorphism types, we show how to produce an infinite family of pairwise nondiffeomorphic nonsymplectic\/ $4$-manifolds belonging to it.  In particular, we prove that there are infinitely many exotic irreducible nonsymplectic smooth structures on $\CP\#3\CPb$, $3\CP\#5\CPb$ and $3\CP\#7\CPb$.
\end{abstract}

\maketitle

\section{Introduction}
\label{sec:intro}

Over the last decade constructing exotic smooth structures on closed oriented simply-connected 4-manifolds with small Euler characteristics has become one of the fundamental problems in differential topology. Here, ``exotic'' means that these manifolds are homeomorphic but not diffeomorphic to standard manifolds, whereas the ``small'' examples we focus on are the ones with either $b_2^+ =1$ or $3$. The most recent history can be split into two periods. The first period was opened by J. Park's paper \cite{Pa:2005} which popularized the rational blowdown technique of Fintushel and Stern \cite{FS: blowdown}, and several constructions of small exotic manifolds relied on an artful use of rational blowdown techniques combined with improved knot surgery tricks \cite{SS, FS: double node, PSS}. More recently, Akhmedov's construction in \cite{A1} triggered the hope that using building blocks with nontrivial fundamental groups could succeed in obtaining exotica on simply-connected $4$-manifolds. These techniques were initially espoused by Fintushel and Stern in \cite{FS: same SW} and later discussed in \cite{Stern} and in \cite{FPS}. The common theme in the recent constructions (\cite{A1, AP, BK:1-3, FPS}) that will also be discussed herein is the manipulations to kill the fundamental group.

The first goal of our paper is to outline a general recipe to obtain small exotic symplectic $4$-manifolds and to fit \emph{all}\/ the recent constructions in \cite{A1,AP,BK:1-3} in this construction scheme. In particular, we hope to make it apparent that seemingly different examples are closely related through a sequence of Luttinger surgeries. The second goal is to calculate the basic classes and the Seiberg-Witten (SW) invariants of these small $4$-manifolds. Using these calculations we show how to obtain infinite families of pairwise nondiffeomorphic manifolds in the homeomorphism type of $\CP\#(2k+1)\CPb$, for $k = 1,\dots,4$, or of $3\CP\# (2l+3)\CPb$, for $l =1,\dots,6$, respectively. We distinguish the diffeomorphism types of these 4-manifolds by comparing their SW invariants. Each of our families will have exactly one symplectic member.

The nature of our constructions differ from the earlier constructions that utilized rational blowdown and knot surgery techniques in the sense that the latter essentially start with elliptic surfaces $E(1)=\CP\#9\CPb$ and $E(2)$ (which has the same Euler characteristic and signature as $3 \CP \# 19 \CPb$), and obtain smaller manifolds from them, whereas our method uses very small manifolds with nontrivial fundamental groups to obtain `bigger' manifolds within the same range (namely $\CP\#(2k+1)\CPb$, for $k = 1,\dots,4$, and $3\CP\# (2l+3)\CPb$, for $l =1,\dots,8$). This aspect also agrees with the  terminology  ``reverse engineering'' suggested for this approach in a general manner by Fintushel and Stern \cite{Stern, FPS}. A problem that we do not deal with here but find interesting is the comparison of examples obtained through these two approaches that run in opposite directions. We hope that the simplicity of our examples and calculations will help with this task.

\medskip
\section{Preliminaries}

\subsection{Symplectic fiber sum and irreducibility}

Let $\Sigma_g$ denote a closed Riemann surface of genus $g > 0$. Since the universal cover of $\Sigma_g$ is contractible, $\Sigma_g$ is acyclic. It follows from the long exact homotopy sequence of a fibration that any $\Sigma_g$ bundle over $\Sigma_h$ is acyclic. In particular, $\pi_2(\Sigma_g\times\Sigma_h)=0$ and hence $\Sigma_g \times \Sigma_h$ is minimal. One new ingredient in our constructions is the following theorem of Michael Usher:

\begin{thm}[Usher \cite{Us}]\label{thm:usher}
Let\/ $X=Y\#_{\Sigma=\Sigma'}Y'$ be the symplectic sum, where the genus $g$ of\/ $\Sigma$ and\/ $\Sigma'$ is strictly positive.
\begin{itemize}
\item[(i)]  If either\/ $Y\setminus\Sigma$ or $Y'\setminus\Sigma'$ contains an embedded symplectic sphere of square $-1$, then $X$ is not minimal.
\item[(ii)]  If one of the summands, say $Y$ for definiteness, admits the structure of an $S^2$-bundle over a surface of genus $g$ such that $\Sigma$ is a section of this $S^2$-bundle, then $X$ is minimal if and only if $Y'$ is minimal.
\item[(iii)] In all other cases, $X$ is minimal.
\end{itemize}
\end{thm}

Combining the above theorem with the following result gives a way to show that simply-connected
symplectic sums are irreducible.

\begin{thm}[Hamilton and Kotschick \cite{HK}]\label{thm:HK}
Minimal symplectic\/ $4$-manifolds with residually finite fundamental groups are irreducible.
\end{thm}

\medskip
\subsection{Seiberg-Witten invariants}
\label{subsec:SW}

In this section we review the basics of Seiberg-Witten invariant (cf.\ \cite{Wi}). The Seiberg-Witten invariant of a smooth closed oriented $4$-manifold $X$\/ is an integer valued function which is defined on the set of Spin$^{c}$ structures on $X$. If we assume that $H_1(X;\Z)$ has no 2-torsion, then there is a one-to-one correspondence between the set of Spin$^{c}$ structures on $X$\/ and the set of characteristic elements of $H^2(X;\Z)$ as follows: To each Spin$^c$ structure $\mathfrak{s}$ on $X$\/ corresponds a bundle of positive spinors $W^+_{\mathfrak{s}}$ over $X$. Let $c(\mathfrak{s}) = c_1(W^+_{\mathfrak{s}}) \in H^2(X;\Z)$. Then each $c(\mathfrak{s})$ is a characteristic element of $H^2(X;\Z)$; i.e. $c_1(W^+_{\mathfrak{s}})$ reduces mod~2 to $w_2(X)$.

In this setup we can view the \emph{Seiberg-Witten invariant}\/ as an integer valued function
\begin{equation*}
\SW_X: \{ k\in H_2(X ; \Z) \mid {\rm PD}(k) \equiv w_2(X) \pmod2 \}
\longrightarrow \Z ,
\end{equation*}
where ${\rm PD}(k)$ denotes the Poincar\'{e} dual of $k$.  The Seiberg-Witten invariant $\SW_X$ is a diffeomorphism invariant when $b_2 ^+(X)>1$ or when $b_2^+(X)=1$ and $b_2^-(X)\leq 9$ (see \cite{Sz: small perturbation} for the $b_2^+=1$ case).  Its overall sign depends on our choice of an orientation of
\begin{equation*}
H^0(X; \R)\otimes\det H_+^2(X;\R)\otimes \det H^1(X;\R).
\end{equation*}

If $\SW_X(\beta)\neq 0$, then we call $\beta$\/ and its Poincar\'{e} dual ${\rm PD}(\beta)\in H^2(X;\Z)$ a {\it basic class}\/ of $X$. It was shown in \cite{Ta} that the canonical class $K_X=-c_1(X, \omega)$ of a symplectic $4$-manifold $(X, \omega)$\/ is a basic class when $b_2^+(X)>1$ with $\SW_{X}(K_X)=1$. It can be shown that, if $\beta$ is a basic class, then so is $-\beta$ with
\begin{equation*}
\SW_X(-\beta)=(-1)^{(e(X) +\sigma(X))/4}\,\SW_X(\beta) ,
\end{equation*}
where $e(X)$ is the Euler characteristic and $\sigma(X)$ is the signature of $X$.  We say that $X$\/ is of \emph{simple type}\/ if every basic class $\beta$\/ of $X$\/ satisfies
\begin{equation*}
\beta^2 = 2e(X) + 3\sigma(X).
\end{equation*}
It was shown in \cite{Ta2} that symplectic 4-manifolds with $b_2^+>1$ are of simple type.
Let $\Sigma\subset X$\/ be an embedded surface of genus $g(\Sigma)>0$.  If $X$\/ is of simple type and $\beta$\/ is a basic class of $X$, we have the following (generalized) \emph{adjunction inequality}\/ (cf.\ \cite{OS:adj ineq}):
\begin{equation*}
[\Sigma]^2 + |\beta\cdot[\Sigma] | \leq 2g(\Sigma)-2.
\end{equation*}

\medskip
\subsection{Surgery on nullhomologus tori and SW invariants}
\label{subsec: FS infinite}

Let $\Lambda$ be a torus of self-intersetion zero inside a $4$-manifold $X$. Choose a framing of the tubular neighborhood $\nu\Lambda$ of $\Lambda$ in $X$, i.e. a diffeomorphism $\nu\Lambda\cong T^2\times D^2$. Given a simple loop $\lambda$ on $\Lambda$, let $S^1_{\lambda}$ be a loop on the boundary $\partial(\nu\Lambda)\cong T^3$ that is parallel to $\lambda$ under the chosen framing.  Let $\mu_{\Lambda}$ denote a meridian circle to $\Lambda$ in $\partial(\nu\Lambda)$. By the $p/q$\/ surgery on $\Lambda$ with respect to $\lambda$, or more simply by a \emph{$(\Lambda,\lambda,p/q)$ surgery}, we mean the closed 4-manifold
\begin{equation*}
X_{\Lambda,\lambda}(p/q) = (X\setminus\nu\Lambda)
\cup_{\varphi} (T^2\times D^2),
\end{equation*}
where the gluing diffeomorphism $\varphi:T^2\times\partial D^2\rightarrow
\partial(X\setminus\nu\Lambda)$ satisifies
\begin{equation*}
\varphi_{\ast}([\partial D^2])= p[\mu_{\Lambda}] + q[S^1_{\lambda}]
\in H_1(\partial(X\setminus\nu\Lambda);\Z).
\end{equation*}
By Seifert-Van Kampen theorem, one easily concludes that
\begin{equation*}
\pi_1(X_{\Lambda,\lambda}(p/q)) = \pi_1(X\setminus\nu\Lambda)/
\langle [\mu_{\Lambda}]^p [S^1_{\lambda}]^q =1\rangle.
\end{equation*}
In what follows, we will be frequently using the following theorem.

\begin{thm}[Fintushel, Park, Stern \cite{FPS}]\label{thm: FS infinite}
Let $X$ be a closed oriented smooth\/ $4$-manifold which contains a nullhomologous torus\/ $\Lambda$ with $\lambda$ a simple loop on\/ $\Lambda$ such that $S^1_{\lambda}$ is nullhomologous in $X\setminus\nu\Lambda$.  If\/
$X_{\Lambda,\lambda}(0)$ has nontrivial Seiberg-Witten invariant, then
the set\/
\begin{equation}\label{set:1/n surgery}
\{ X_{\Lambda,\lambda}(1/n) \mid n = 1,2,3,\dots \}
\end{equation}
contains infinitely many pairwise nondiffeomorphic\/ $4$-manifolds. Furthermore, if $X_{\Lambda,\lambda}(0)$ only has just one Seiberg-Witten basic class up to sign, then every pair of\/ $4$-manifolds in\/ $(\ref{set:1/n surgery})$ are nondiffeomorphic.
\end{thm}

\begin{rmk} \label{rmk: two classes}
Note that $X=X_{\Lambda,\lambda}(1/0)$.  Let $T$\/ be the core torus of the 0 surgery $X_{\Lambda,\lambda}(0)$.  If $k_0$ is a characteristic element
of $H_2(X_{\Lambda,\lambda}(0);\Z)$ satisfying $k_0\cdot[T]=0$, then $k_0$ gives rise to unique characteristic elements $k\in H_2(X;\Z)$ and $k_n\in H_2(X_{\Lambda,\lambda}(1/n);\Z)$.  The product formula in \cite{MMS} then gives
\begin{equation}\label{eq:MMS}
\SW_{X_{\Lambda,\lambda}(1/n)}(k_n) = \SW_X (k) +
n\sum_{i\in \Z}\SW_{X_{\Lambda,\lambda}(0)}(k_0 + 2i[T]).
\end{equation}
Let us now assume that $X_{\Lambda,\lambda}(0)$ has only one basic class up to sign and this basic class is not a multiple of $[T]$.  Under these assumptions, the infinite sum in (\ref{eq:MMS}) only contains at most one nonzero summand.  If we further assume that $X$\/ and $X_{\Lambda,\lambda}(0)$ are both symplectic, then the adjunction inequality implies that the only basic class of $X$\/ and $X_{\Lambda,\lambda}(0)$ is the canonical class up to sign.  Under all these assumptions, it follows that $X_{\Lambda,\lambda}(1/n)$ also has only one basic class up to sign for every $n\geq 1$.
\end{rmk}

\medskip
\subsection{Luttinger surgery}
\label{subsec:Luttinger}

Luttinger surgery is a special case of $p/q$\/ surgery on a self-intersection zero torus $\Lambda$ described in the previous subsection. It was first studied in \cite{Lu} and then in \cite{ADK} in a more general setting.
Assume that $X$\/ is a symplectic 4-manifold with a symplectic form $\omega$, and that the torus $\Lambda$ is a Lagrangian submanifold of $X$. It is well-known that there is a canonical framing of $\nu\Sigma\cong T^2\times D^2$, called the
{\it Lagrangian framing}, such that $T^2\times\{ x \}$ corresponds to a Lagrangian submanifold of $X$\/ for every $x\in D^2$. Given a simple loop $\lambda$ on $\Lambda$, let $S^1_{\lambda}$ be a simple loop on $\partial(\nu\Lambda)$ that is parallel to $\lambda$ under the Lagrangian framing. For any integer $m$, the $(\Lambda,\lambda,1/m)$ {\it Luttinger surgery}\/ on $X$\/ will be $X_{\Lambda,\lambda}(1/m)$, the $1/m$\/ surgery on $\Lambda$ with respect to $\lambda$ and
the Lagrangian framing.  Note that our notation is different from the one in \cite{ADK} wherein $X_{\Lambda,\lambda}(1/m)$ is denoted by $X(\Lambda,\lambda,m)$.  It is  shown in \cite{ADK} that $X_{\Lambda,\lambda}(1/m)$ possesses a symplectic form that restricts to the original symplectic form $\omega$ on $X\setminus\nu\Lambda$.  In this paper, we will \emph{only}\/ look at Luttinger surgeries when $m=\pm 1=1/m$, so there should be no confusion in notation.

\begin{rmk}
In Section~\ref{sec:infinite families} and the Appendix, we will also be looking at a \emph{non}-Luttinger $(\Lambda,\lambda,-n)$ surgery $X_{\Lambda,\lambda}(-n)$ for a Lagrangian torus $\Lambda$ equipped with the Lagrangian framing and a positive integer $n\geq 2$.
\end{rmk}

\medskip
\section{The general construction scheme}
\label{sec:general}

Here we outline a general construction scheme to construct simply-connected minimal symplectic $4$-manifolds with small Euler characteristics. This is an incidence of the ``reverse engineering'' (\cite{Stern, FPS}) idea applied to certain symplectic manifolds. Any example using this scheme and homeomorphic to $\CP \# n \CPb$ (for $n > 0$) and $m \CP \# n \CPb$ (for $m >0$) can be distinguished from the latter standard manifolds by comparing their symplectic structures or their Seiberg-Witten invariants, respectively. Recall that, $\CP \# n \CPb$ (for $n > 0$) are nonminimal, and  $m \CP \# n \CPb$ (for $m > 0$) all have vanishing Seiberg-Witten invariants, unlike the minimal symplectic 4-manifolds that we produce. Our approach will allow us to argue easily how all 4-manifolds obtained earlier in \cite{A1, AP, BK:1-3} arise from this construction scheme, and in particular we show how seemingly different examples rely on the very same idea.

The only building blocks we need are the products of two Riemann surfaces. In fact, it suffices to consider multiple copies of $S^2\times T^2$ and $T^2 \times T^2$, since all the other product manifolds except for $S^2 \x S^2$ (which we will not use here) can be obtained by fiber summing copies of these manifolds appropriately. Note that any such manifold is a minimal symplectic manifold. Both $S^2\times T^2$ and $T^4 = T^2 \times T^2$ can be equipped with product symplectic forms where each factor is a symplectic submanifold with self-intersection zero. Denote the standard generators of $\pi_1(T^4)$ by $a$, $b$, $c$\/ and $d$, so that $H_2(T^4;\Z) \cong \Z^6$ is generated by the homology classes of two symplectic tori $a \times b$\/ and $c \times d$, and four Lagrangian tori $a\times c$, $a \times d$, $b \times c$\/ and $b \times d$\/ with respect to the product symplectic form on $T^4$ that we have chosen. The intersection form splits into three hyperbolic pairs: $a \times b$\/
and $c \times d$, $a \times c$\/ and $b \times d$, $a \times d$\/ and $b \times c$. Finally, note that all four Lagrangian tori can be pushed off to nearby Lagrangians in its standard Weinstein neighborhoods so that they lie in the complement of small tubular neighborhoods of the two symplectic tori. With a little abuse of notation (which will be remembered in our later calculations of fundamental groups), we will still denote these parallel Lagrangian tori with the same letters.

In order to produce an exotic copy of a target manifold $Z$, we first perform blow-ups and symplectic fiber sums to obtain an intermediate manifold $X'$. Whenever a piece is blown-up, we make sure to fiber sum that piece along a symplectic surface that intersects each exceptional sphere positively at one point. This allows us to employ Theorem~\ref{thm:usher} to conclude that $X'$ is minimal. We want this intermediate manifold to satisfy the following two properties:
\begin{itemize}
\item[(I)] $X'$ should have the same signature and Euler characteristic as $Z$.
\item[(II)] If $r$\/ is the rank of the maximal subspace of $H_2(X';\Z)$ generated by
homologically essential Lagrangian tori, then we should have $r \geq s = 2 b_1(X')= b_2(X') - b_2(Z)$.
\end{itemize}
Moreover, we generally desire to have $\pi_1(X') = H_1(X';\Z)$ for the reasons that will become apparent below. However, surprisingly one can also handle some examples where $\pi_1(X')$ is not abelian. (See for example Section \ref{subsec:3-5}, or \cite{FPS}.)

Finally, we carefully perform $s/2$\/ Luttinger surgeries to kill $\pi_1(X')$ and obtain a simply-connected symplectic $4$-manifold $X$. Note that these surgeries can easily be chosen to obtain a manifold with $b_1 = 0$. However, determining the correct choice of Luttinger surgeries in this last step to kill the fundamental group
completely is a much more subtle problem. This last part is certainly the hardest part of our approach, at least for the
`smaller' constructions. The reader might want to compare below the complexity of our fundamental group calculations for $\CP\#(2k+1)\CPb$, for $k = 1,\dots,4$ as $k$ gets smaller.

In order to compute and effectively kill the fundamental group of the resulting manifold $X$, we will do the Luttinger surgeries in our building blocks as opposed to doing them in $X'$. This is doable, since the Lagrangian tori along which we perform Luttinger surgeries lie away from the symplectic surfaces that are used in any symplectic sum constructions, as well as the blow-up regions. In other words, one can change the order of these operations while paying extra attention to the $\pi_1$ identifications. Having the $\pi_1$ calculations of the pieces in hand, we can use Seifert-Van Kampen theorem repeatedly to calculate the fundamental group of our exotic candidate $X$.

Below, we will work out some concrete examples, where we construct minimal symplectic 4-manifolds homeomorphic to $\CP\#(2k+1)\CPb$, for $k = 1, \dots, 4$, and $3 \CP \# (2l+3) \CPb$, for $l =1,\dots, 6$. We hope that the reader will have a better understanding of the recipe we have given here by looking at these examples. Another essential observation that is repeatedly used in our arguments below is the interpretation of some manifold pieces used in \cite{A1,AP} as coming from Luttinger surgeries on $T^4$, together with the description of their fundamental groups.
This is proved in the Appendix.  A concise history of earlier constructions will be given at the beginning of each subsection.

\medskip
\subsection{A new description of an exotic $\CP\# 9\CPb$}
\label{subsec:1-9}

The first example of an exotic smooth structure on the elliptic surface $E(1)=\CP\#9\CPb$, and in fact the first exotic smooth structure on any closed topological 4-manifold, was constructed by Donaldson in \cite{Do}. Donaldson's example was the Dolgachev surface $E(1)_{2,3}$. Later on, Friedman showed that $\{E(1)_{p,q}\mid \gcd (p,q) =1\}$
contains infinitely many nondiffeomorphic 4-manifolds (cf.\ \cite{Fri}). In \cite{FS: knots links} Fintushel and Stern have shown that knot surgered manifolds $E(1)_K$ give infinitely many irreducible smooth structures on $E(1)=\CP\#9\CPb$.

Consider $S^2 \times T^2 = S^2 \times (S^1 \times S^1)$ equipped with its product symplectic form, and denote the last two circle factors by $x$\/ and $y$. One can take the union of three symplectic surfaces $(\{s_1\} \times T^2) \cup (S^2 \times \{t\} )\cup (\{s_2\} \times T^2)$ in $S^2 \times T^2$, and resolve the two double points symplectically.  This yields a genus two symplectic surface in $S^2\times T^2$ with self-intersection four. Symplectically blowing up
$S^2 \times T^2$ along these four intersection points and taking the proper transform, we obtain a symplectic genus two surface $\Sigma$ in $Y= (S^2 \times T^2) \# 4 \CPb$.  Note that the inclusion induced homomorphism from
$\pi_1(\Sigma) = \langle a, b, c, d \mid [a,b][c,d]=1 \rangle$ into $\pi_1(Y) = \langle x, y \mid
[x,y]=1 \rangle$ maps the generators as follows:
\begin{equation*}
a \mapsto x, \; b \mapsto y,\; c \mapsto x^{-1},\; d \mapsto y^{-1} .
\end{equation*}

Let us run the same steps in a second copy of $S^2 \times T^2$ and label every object with a prime symbol at the end.  That is, $Y' = (S^2 \times T^2)\# 4 \CPb$, $\Sigma'$ is the same symplectic genus two surface described above with $\pi_1$ generators $a', b', c', d'$, and finally let $x', y'$ denote the generators of the $\pi_1(Y')$. Let $X$\/ be the symplectic fiber sum of $Y$\/ and $Y'$ along $\Sigma$ and $\Sigma'$ via a diffeomorphism that extends the
orientation-preserving diffeomorphism $\phi: \Sigma \to \Sigma'$, described by:
\begin{equation*}
a \mapsto a'b',\; b \mapsto (a')^{-1},\; c \mapsto c',\; d \mapsto d' .
\end{equation*}

The Euler characteristic of $X$\/ can be computed as $e(X) = 4 + 4 - 2(2- 2\cdot2) = 12$, and the Novikov additivity gives the signature $\sigma(X)= -4 + (-4) = -8$, which are exactly the Euler characteristic
and the signature of $Z=\CP\#9\CPb$. We claim that $X$\/ is already simply-connected and thus no Luttinger surgery is needed. Note that $\pi_1(Y\setminus\nu\Sigma)=\pi_1(Y)$ since a meridian circle of $\Sigma$ bounds a punctured exceptional sphere from one of the four blowups. Using Seifert-Van Kampen theorem, we see that
\begin{eqnarray*}
\pi_1(X)&=& \langle x,y,x',y' \mid [x,y]=[x',y']=1,\\
&& x = x'y',\, y = (x')^{-1},\, x^{-1} = (x')^{-1},\, y^{-1}= (y')^{-1}
\rangle.
\end{eqnarray*}
We conclude that $x = x'$, $y = y'$, $y = x^{-1}$. Thus $x = x'y'$ implies $y =1$, and in turn $x=1$. So $\pi_1(X) = 1$. Hence by Freedman's classification theorem for simply-connected topological 4-manifolds (cf.\ \cite{Fre}), $X$\/ is homeomorphic to $E(1)$. However, $X$\/ is irreducible by Theorem~\ref{thm:HK}, and therefore $X$\/ is not diffeomorphic to $E(1)$. The 4-manifold $X$\/ we obtained here can be shown to be the knot surgered manifold $E(1)_K$, where the knot $K$\/ is the trefoil
(cf.\ \cite{FS: same SW}).

Alternatively we could construct the above manifold in the following way.  First we symplectically sum two copies of $(S^2 \times T^2) \# 4 \CPb$ along $\Sigma$ and $\Sigma'$ via a map that directly identifies the generators $a, b, c, d$\/ with $a', b', c', d'$
in that order. Call this symplectic 4-manifold $X'$ and observe that while the characteristic numbers $e$ and $\sigma$ are the same as above, this manifold has $\pi_1(X') = H_1(X';\Z) \cong \Z^2$ and $H_2(X';\Z)$ has four additional classes that do not occur in $X$. These classes are as follows.  Inside $((S^2 \times T^2) \# 4 \CPb) \setminus \nu \Sigma$, there are cylinders $C_a$ and $C_b$ with
\begin{eqnarray*}
\partial C_a &=& a \cup c, \quad \partial C_b = b \cup d.
\end{eqnarray*}
Similiarly we obtain cylinders $C'_a$ and $C'_b$ in the second copy of $((S^2 \times T^2) \# 4 \CPb) \setminus \nu \Sigma'$. Thus we can form the following internal sums in $X'$:
\begin{eqnarray*}
\Sigma_a &=&  C_a \cup C_a' ,\quad  \Sigma_b = C_b \cup C_b' .
\end{eqnarray*}
These are all tori of self-intersection zero. Let $\mu$ denote a meridian of $\Sigma$, and let $R_{a}=\tilde{a} \times \mu$, and $R_b = \tilde{b} \times \mu$\/ be the `rim tori', where $\tilde{a}$\/ and $\tilde{b}$\/ are suitable parallel copies of the generators $a$\/ and $b$. Note that $[R_a]^2=[R_b]^2=[\Sigma_a]^2= [\Sigma_b]^2=0$, and $[R_a]\cdot[\Sigma_b]=1=[R_b]\cdot[\Sigma_a]$.

Observe that these rim tori are in fact Lagrangian.  One can show that the effect of two Luttinger surgeries $(R_a, \tilde{a} , -1)$ and $(R_b, \tilde{b}, -1)$ is the same as changing the gluing map that we have used in the symplectic sum to the gluing map $\phi$\/ in the first construction. This second viewpoint is the one that will fit in with our construction of an infinite family of pairwise nondiffeomorphic smooth structures in Section~\ref{sec:infinite families}.

\medskip
\subsection{A new construction of an exotic $\CP\#7\CPb$}
\label{subsec:1-7}

The first example of an exotic $\CP\#7\CPb$ was constructed by J. Park in \cite{Pa:2005} by using rational blowdown (cf.\ \cite{FS: blowdown}), and the Seiberg-Witten invariant calculation in \cite{OS} shows that it is irreducible. Infinitely many exotic examples were later constructed by Fintushel and Stern in \cite{FS: double node}. All of their constructions use the rational blowdown technique. Here we construct another irreducible symplectic $4$-manifold homemorphic but not diffeomorphic to $\CP\#7\CPb$ using our scheme, and thus without using any rational blowdown.

We equip $T^4 = T^2 \times T^2$ and $S^2 \times T^2$ with their product symplectic forms. The two orthogonal symplectic tori in $T^4$ can be used to obtain a symplectic surface of genus two with self-intersection two. Symplectically blowing-up at these self-intersection points we obtain a new symplectic surface $\Sigma$ of genus two with trivial normal bundle in $Y = T^4 \# 2 \CPb$. The generators of $\pi_1(T^4 \# 2 \CPb)$ are the circles $a, b, c, d$, and the inclusion induced homomorphism from $\pi_1(\Sigma)$ to
\begin{equation*}
\pi_1(Y)= \langle a, b, c, d \mid
[a,b]=[a,c]=[a,d]=[b,c]=[b,d]=[c,d]=1 \rangle
\end{equation*}
is surjective. Indeed the four generators of $\pi_1(\Sigma)$ are mapped onto $a, b, c, d$\/ in $\pi_1(Y)$, respectively.

On the other hand, as in Subsection~\ref{subsec:1-9}, we can start with $S^2\times T^2$ and get a symplectic genus two surface $\Sigma'$ in $Y'=(S^2\times T^2)\# 4\CPb$. Once again $\pi_1(Y') = \langle x, y \mid [x,y] =1 \rangle$ and the generators $a',b',c',d'$ of $\pi_1(\Sigma')$ are identified with $x,y,x^{-1},y^{-1}$, respectively.

We take the symplectic sum of $Y$\/ and $Y'$ along $\Sigma$ and $\Sigma'$ given by a diffeomorphism that extends the
identity map sending $a \mapsto a', b \mapsto b', c \mapsto c', d \mapsto d'$ to obtain an intermediate 4-manifold $X'$. The Euler characteristic can be computed as $e(X') = 2 + 4 + 4 = 10$, and the Novikov additivity gives $\sigma(X')= -2 + (-4) = -6$, which are the characteristic numbers of $\CP\#7\CPb$. Since exceptional spheres intersect $\Sigma$ and $\Sigma'$ transversally once, we have $\pi_1(Y\setminus\nu\Sigma) \cong \pi_1(Y)$ and $\pi_1(Y'\setminus\nu\Sigma') \cong \pi_1(Y')$. Using Seifert-Van Kampen theorem, we compute that
\begin{eqnarray*}
\pi_1(X') &=& \langle a,b,c,d,x,y \mid [a,b]=[a,c]=[a,d]=[b,c]=[b,d]=[c,d]=1, \\
&& [x,y]=1,\, a=x,\, b=y,\, c=x^{-1},\, d=y^{-1} \rangle.
\end{eqnarray*}
Thus $\pi_1(X') = \langle x, y \mid [x,y]=1 \rangle\cong\Z^2$, and it follows that $b_2(X') = 12$ from our Euler characteristic calculation above. The four homologically essential Lagrangian tori in $T^4$ are also contained in $X'$, and thus one can see that condition (II) is satisfied.

The two Luttinger surgeries we choose are $-1$ surgery on $\tilde{a} \times c$\/ along $\tilde{a}$\/ and another $-1$ surgery on $\tilde{b} \times c$\/ along $\tilde{b}$.   Here,
$\tilde{a}$\/ and $\tilde{b}$\/ are suitable parallel copies of the generators $a$\/ and $b$, respectively.  We claim that the manifold $X$\/ we obtain after these two Luttinger surgeries is simply-connected.
To prove our claim, we observe that these two Luttinger surgeries could be first made in the $T^4$ piece that we had at the very beginning.  This is because both Lagrangian tori $\tilde{a} \times c$\/ and $\tilde{b} \times c$\/ lie in the complement of $\Sigma$.
By our observation in the Appendix, the result of these two Luttinger surgeries in $T^4$ is diffeomorphic to $S^1 \times M_K$. Observe that $\pi_1((S^1 \times M_K)\# 2\CPb
\setminus \nu\Sigma) \cong \pi_1(S^1 \times M_K )$, which is (cf.\ \cite{AP} and (\ref{pi_1 relations from two surgeries})--(\ref{monodromy2}) in the Appendix)
\begin{eqnarray*}
\langle a,b,c,d \mid [a,b]=[c,a]=[c,b]=[c,d]=1,\,
d a d^{-1}=[d,b^{-1}], \, b=[a^{-1},d] \rangle.
\end{eqnarray*}
As before, $\pi_1(((S^2 \times T^2)\# 4\CPb) \setminus \nu\Sigma')
\cong \pi_1(S^2 \times T^2) = \langle x, y \mid [x,y]=1 \rangle$.
Therefore by Seifert-Van Kampen theorem,
\begin{eqnarray*}
\pi_1(X) &=& \langle a,b,c,d,x,y \mid [a,b]=[c,a]=[c,b]=[c,d]=1,\\
&&d a d^{-1}=[d,b^{-1}], \, b=[a^{-1},d],\, [x,y]=1, \\
&& a=x,\, b=y,\, c=x^{-1},\, d=y^{-1} \rangle.
\end{eqnarray*}

Thus $x$\/ and $y$\/ generate the whole group, and by direct substitution we see that $y^{-1} x y = [y^{-1},y^{-1}]=1$ and $y = [x^{-1},y^{-1}]$. The former gives $x = 1$, and the latter then yields $y =1$. Hence $\pi_1(X) = 1$. Therefore by Freedman's classification theorem (cf.\ \cite{Fre}), $X$\/ is homeomorphic to $\CP\#7\CPb$.  Since the latter is not irreducible, $X$\/ is an exotic copy of it.

\medskip
\subsection{A new construction of an exotic $\CP\#5\CPb$}
\label{subsec:1-5}

The first example of an exotic $\CP\# 5\CPb$ was obtained by J. Park, Stipsicz and Szab\'{o} in \cite{PSS},
combining the double node neighborhood surgery technique discovered by Fintushel and Stern (cf.\ \cite{FS: double node}) with rational blowdown.  Fintushel and Stern also constructed similar examples using the same techniques
in \cite{FS: double node}. The first exotic symplectic $\CP\# 5\CPb$ was constructed by the first author in \cite{A1}. Here, we present another construction with a much simpler $\pi_1$ calculation,
using our construction scheme.

As in Subsection~\ref{subsec:1-7}, we construct a symplectic surface $\Sigma$ of genus two with trivial normal bundle in $Y=T^4 \# 2 \CPb$. Let us use the same notation for the fundamental groups as above. Take another copy
$Y'=T^4 \# 2 \CPb$, and denote the same genus two surface by $\Sigma'$, while using the prime notation for all corresponding fundamental group elements.

We obtain a new manifold $X'$ by taking the symplectic sum of $Y$\/ and $Y'$ along $\Sigma$ and $\Sigma'$ determined by the map $\phi: \Sigma \to \Sigma'$ that satisfies:
\begin{equation}\label{eq:phi images}
a \mapsto c',\; b \mapsto
d',\; c \mapsto a',\;  d \mapsto b' .
\end{equation}
By Seifert-Van Kampen theorem, one can easily verify that $\pi_1(X') \cong \Z^4$ generated by, say $a,
b, a', b'$. The characteristic numbers we get are:  $e(X') = 2+ 2 + 4 = 8$ and $\sigma(X') = -2 +(-2) = -4$, the characteristic numbers of $\CP\# 5\CPb$. Finally the homologically essential Lagrangian tori in the initial $T^4$ copies can be seen to be contained in $X'$ with the same properties. Thus $r \geq 8 = 2 b_1(X') = b_2(X') - b_2(\CP\# 5\CPb)$, so our condition (II) is satisfied.

We perform the following four Luttinger surgeries on pairwise disjoint Lagrangian tori:
\begin{gather*} (\tilde{a} \times c, \tilde{a}, -1), \ \ (\tilde{b} \times c, \tilde{b}, -1) ,\ \ (\tilde{a}' \times c', \tilde{a}', -1) ,\ \  (\tilde{b}' \times c', \tilde{b}', -1).
\end{gather*}

It is quite simple to see that the resulting symplectic 4-manifold $X$\/ satisfies $H_1(X;\Z) = 0$. Using the Appendix again, after changing the order of operations and assuming that we have done the Luttinger surgeries at the very beginning, we can view $X$\/ as the fiber sum of two copies of $(S^1 \times M_K)\# 2\CPb$ along the identical genus two surface $\Sigma$ where the gluing map switches the symplectic bases for $\Sigma$ as in (\ref{eq:phi images}). Thus, using Seifert-Van Kampen's theorem as above, we can see that
\begin{eqnarray*}
\pi_1(X) &=& \langle a,b,c,d,a',b',c',d' \mid
[a,b]=[c,a]=[c,b]=[c,d]=1,\\
&& dad^{-1} =  [d,b^{-1}],\, b=[a^{-1},d] ,\,
[a',b']=[c',a']=[c',b']=[c',d']=1, \\
&& d'a'(d')^{-1} = [d',(b')^{-1}],\, b'=[(a')^{-1},d'],\, \\
&& a = c',\, b = d',\,
c =a',\, d= b' \rangle.
\end{eqnarray*}
Now $b'=[(a')^{-1},d']$ can be rewritten as $d=[c^{-1},b]$.  Since $b$\/ and $c$\/ commute, $d=1$.  The relations $dad^{-1} =  [d,b^{-1}]$ and $b=[a^{-1},d]$ then quickly implies that $a = 1$ and $b = 1$, respectively.
Lastly, $d'a'(d')^{-1} = [d',(b')^{-1}]$ is $bcb^{-1} = [b,d^{-1}]$, so $c=1$ as well. Since $a, b, c, d$\/ generate $\pi_1(X)$, we see that $X$\/ is simply-connected.  By similar arguments as before, $X$\/ is an irreducible symplectic $4$-manifold that is homeomorphic but not diffeomorphic to $\CP\# 5\CPb$.

\medskip
\subsection{An exotic $\CP\#3\CPb$, revisited}
\label{subsec:1-3}

The first exotic irreducible symplectic smooth structure on $\CP\#3\CPb$ was constructed in \cite{AP}. An alternative construction appeared in \cite{BK:1-3}. In light of the Appendix, these two constructions can be seen to differ only by the choice of three out of six Luttinger surgeries.

Let us demonstrate how the construction of an exotic symplectic $\CP \# 3 \CPb$ in \cite{AP} fits into our recipe.  We will use three copies of the $4$-torus, $T^4_1$, $T^4_2$ and $T^4_3$. Symplectically fiber sum the first two along
the 2-tori $a_1 \times b_1$ and $a_2 \times b_2$ of self-intersection zero, with a gluing map that identifies $a_1$ with $a_2$ and $b_1$ with $b_2$. Clearly we get $T^2 \times \Sigma_2$, where the symplectic genus 2 surface $\Sigma_2$ is obtained by gluing together the orthogonal punctured symplectic tori $(c_1 \times d_1)\setminus D^2$ in $T^4_1$ and $(c_2 \times d_2)\setminus D^2$ in $T^4_2$. Here, $\pi_1(T^2 \times \Sigma_2)$ has six generators $a_1 = a_2$, $b_1 = b_2$, $c_1$, $c_2$, $d_1$ and $d_2$ with relations $[a_1,b_1]=1$, $[c_1,d_1][c_2,d_2]=1$ and moreover $a_1$ and $b_1$
commute with all $c_i$ and $d_i$. The two symplectic tori $a_3 \times b_3$ and $c_3 \times d_3$ in $T^4_3$ intersect at one point, which can be smoothened to get a symplectic surface of genus two. Blowing up $T^4_3$ twice at the self-intersection points of this surface as before, we obtain a symplectic genus two surface $\Sigma'$\/ of self-intersection zero.

Next we take the symplectic fiber sum of $Y=T^2 \times \Sigma_2$ and $Y'=T^4_3 \# 2 \CPb$ along the surfaces $\Sigma_2$ and $\Sigma'$, determined by a map that sends the circles $c_1, d_1, c_2, d_2$ to $ a_3, b_3, c_3, d_3$
in the same order. By Seifert-Van Kampen theorem, the fundamental group of the resulting manifold $X'$ can be seen to be generated by $a_1, b_1, c_1, d_1, c_2$ and $d_2$, which all commute with each other.  Thus $\pi_1(X')$ is isomorphic to $\Z^6$. It is easy to check that $e(X') = 6$ and $\sigma(X') = -2$, which are also the characteristic numbers of $\CP \# 3 \CPb$.

Now we perform six Luttinger surgeries on pairwise disjoint Lagrangian tori:
\begin{gather*}
(a_1 \times \tilde{c}_1, \tilde{c}_1, -1), \ \ (a_1 \times \tilde{d}_1, \tilde{d}_1, -1),\ \ (\tilde{a}_1 \times c_2, \tilde{a}_1, -1), \\
(\tilde{b}_1 \times c_2, \tilde{b}_1, -1),\ \ (c_1 \times \tilde{c}_2, \tilde{c}_2, -1),\ \ (c_1 \times \tilde{d}_2, \tilde{d}_2, -1).
\end{gather*}
Afterwards we obtain a symplectic 4-manifold $X$\/ with $\pi_1(X)$ generated by $a_1$, $b_1$, $c_1$, $d_1$, $c_2$, $d_2$ with relations:
\begin{gather*}
[b_1,d_1^{-1}]=b_1c_1b_1^{-1},\ \ [c_1^{-1},b_1]=d_1,\ \ [d_2, b_1^{-1}]=d_2a_1d_2^{-1},\\
[a_1^{-1},d_2]=b_1,\ \ [d_1,d_2^{-1}]=d_1c_2d_1^{-1},\ \ [c_2^{-1},d_1]=d_2,
\end{gather*}
and all other commutators are equal to the identity. Since $[b_1, c_2] = [c_1, c_2] =1$, $d_1= [c_1^{-1},b_1]$ also commutes with $c_2$. Thus $d_2 = 1$, implying $a_1 = b_1 = 1$. The last identity implies $c_1=d_1=1$, which in turn implies $c_2 = 1$.

Hence $X$\/ is simply-connected and since these surgeries do not change the characteristic numbers, we have it homeomorphic to $\CP \# 3 \CPb$.  Since $Y$\/ is minimal and the exceptional spheres in $Y'$ intersect $\Sigma'$, Theorem~\ref{thm:usher} guarantees that $X'$ is minimal.  It follows from Theorem~\ref{thm:HK} that $X$\/ is an irreducible symplectic $4$-manifold which is not diffeomorphic to $\CP\# 3\CPb$.

\medskip
\subsection{New constructions of exotic $3\CP\# (2l+3)\CPb$ for $l = 2,\dots,6$}
\label{subsec:3-(2l+3)}

We begin by proving the following theorem.

\begin{thm} \label{thm: lattice jump}
Let $X$ be a simply-connected minimal symplectic\/ $4$-manifold which is not a sphere bundle over a Riemann surface and such that $X$ contains a genus two symplectic surface of self-intersection zero. Then $X$ can be used to construct simply-connected irreducible symplectic\/ $4$-manifolds $Z'$ and $Z''$ satisfying:
\begin{eqnarray*}
(b_2^+(Z'), b_2^-(Z')) &=& (b_2^+(X)+2, b_2^-(X) +4),\\
(b_2^+(Z''), b_2^-(Z'')) &=& (b_2^+(X)+2, b_2^-(X) +6).
\end{eqnarray*}
\end{thm}

\begin{proof}
Let us denote the genus two  symplectic surface of self-intersection $0$ in $X$\/ by $\Sigma_2$. By our assumptions, the complement $X \setminus \nu \Sigma_2$ does not contain any exceptional spheres. Take $T^4=T^2\times T^2$ equipped with a product symplectic form, with the genus two symplectic surface that is obtained from the two orthogonal symplectic tori after
resolving their singularities. After symplectically blowing up $T^4$ at two points on this surface, we get a symplectic genus two surface $\Sigma_2'$ of self-intersection $0$ in $T^4 \# 2 \CPb$, and it is clear that $(T^4 \# 2 \CPb) \setminus \nu \Sigma_2'$ does not contain any exceptional spheres either. Since we also assumed that $X$\/ was
not a sphere bundle over a Riemann surface, it follows from Theorems~\ref{thm:usher} and \ref{thm:HK} that the $4$-manifold $Z'$ obtained as the symplectic sum of $X$\/ with $T^4 \# 2 \CPb$ along $\Sigma_2$ and $\Sigma_2'$ is minimal and hence irreducible.

Next we take $S^2 \times T^2$ with its product symplectic form, and as before consider the genus two symplectic surface obtained from two parallel copies of the symplectic torus component and a symplectic sphere component, after symplectically resolving their intersections. Symplectically blowing up $S^2 \x T^2$ on four points on this surface, we get a new symplectic genus 2 surface $\Sigma_2''$ with self-intersection $0$ in $(S^2 \x T^2) \# 4 \CPb$. Although this second piece $(S^2 \times T^2) \# 4 \CPb$ is an $S^2$ bundle over a Riemann surface, the surface $\Sigma''_2$ cannot be a section of this bundle. Moreover, it is clear that $((S^2 \times T^2) \# 4 \CPb) \setminus \nu \Sigma''_2$ does not contain any exceptional spheres.
Hence, applying Theorems~\ref{thm:usher} and \ref{thm:HK} again, we see that the $4$-manifold $Z''$ obtained as the
symplectic sum of $X$\/ with $(S^2 \times T^2)\# 4\CPb$ along $\Sigma_2$ and $\Sigma_2''$ is minimal and irreducible.

It is a straightforward calculation to see that $(e(Z'), \sigma(Z')) = (e(X) + 6, \sigma(X) - 2)$ and $(e(Z''), \sigma(Z'')) = (e(X) + 8, \sigma(X) - 4)$. Note that the new meridian in $X \setminus \nu \Sigma_2$ dies after the fiber sum since the meridian of $\Sigma_2'$ in $T^4 \# 2\CPb$ can be killed along any one of the two exceptional spheres.  The same argument works for the fiber sum with $(S^2 \x T^2) \# 4 \CPb$. Hence Seifert-Van Kampen's theorem implies that $\pi_1(Z') = \pi_1(Z'') = 1$. Our claims about $b^+_2$ and $b^-_2$ follow immediately.
\end{proof}

\begin{cor} \label{cor: 3CPs}
There are exotic\/ $3 \CP\# (2l+3) \CPb$, for $l = 2, \dots, 6$, which are all irreducible and symplectic.
\end{cor}

\begin{proof}
We observe that each one of the irreducible symplectic $\CP \# (2k+1) \CPb$ ($k = 1, \ldots, 4$) we obtained above contains at least one symplectic genus two surface of self-intersection zero. (Also see Section~\ref{sec:infinite families} for more detailed description of these surfaces.) To be precise, let us consider the genus two surface
$\Sigma$ which is a parallel copy of the genus two surface used in the last symplectic sum in any one of our constructions. Since these exotic $4$-manifolds are all minimal, they cannot be the total space of a sphere bundle over a Riemann surface with any blow-ups in the fibers.  Also they cannot be homeomorphic to either $F \times S^2$ or $F \tilde{\times} S^2$ for some Riemann surface $F$, because of their intersection forms. Therefore we see that assumptions of Theorem~\ref{thm: lattice jump} hold. It quickly follows from Theorem~\ref{thm: lattice jump} that we can obtain irreducible symplectic $4$-manifolds homeomorphic to $3\CP \# (2l + 3) \CPb$, for $l = 2, \dots, 6$.
\end{proof}

\begin{rmk} \label{the rest}
Using the generic torus fiber and a sphere section of self-intersection $-1$ in an elliptic fibration on $E(1)=\CP\#9\CPb$, one can form a smooth symplectic torus $T_1$ of self-intersection $+1$ in $E(1)$. As each one of our exotic $\CP \# (2k+1)\CPb$ for $k = 1, \ldots, 4$ contains at least one symplectic torus of self-intersection $-1$ (these tori are explicitly described in Section~\ref{sec:infinite families}), we can symplectically fiber sum each exotic $\CP\#(2k+1)\CPb$ with $E(1)$ along a chosen torus of self-intersection $-1$ and $T_1$ to obtain irreducible symplectic 4-manifolds that are homeomorphic but not diffeomorphic to $3 \CP \# (2k + 11) \CPb$ for $k = 1, \dots, 4$.  Also compare with
the second paragraph following Question 1 in Section~\ref{sec:questions}.
\end{rmk}

\medskip
\subsection{Construction of an exotic $3\CP\# 5\CPb$, revisited}
\label{subsec:3-5}

This time we use the pieces $\Sigma_2 \times \Sigma_2$ and $T^4$. We view $\Sigma_2 \times \Sigma_2$ as a result of the obvious symplectic fiber sums of four copies of $T^4$, and thus it contains eight pairs of homologically essential Lagrangian tori, each copy of $T^4$ contributing two pairs. Since the details are very much the same as in the previous subsections, we leave it to the reader to verify that the consruction of exotic $3\CP\# 5\CPb$ in \cite{AP} is indeed obtained from our general construction scheme applied to the intermediate symplectic 4-manifold $X'$ which is the fiber sum of $\Sigma_2 \times \Sigma_2$ with $T^4 \# 2 \CPb$ along one of the symplectic components $\Sigma$ in the first piece and the usual genus two surface $\Sigma'$ in $T^4\#2\CPb$.  One then proceeds with four pairs of Luttinger surgeries along eight of the aforementioned Lagrangian tori in the $\Sigma_2\times\Sigma_2$ half.  Each pair of Luttinger surgeries is to be performed in a fixed copy of $T^4$ as described in the Appendix. Notice that initially $X'$ does not have an abelian fundamental group.

\medskip
\section{Infinite families of small exotic 4-manifolds}
\label{sec:infinite families}

In this section we will show how to construct infinite families of pairwise nondiffeomorphic 4-manifolds that are
homeomorphic to one of $\CP\#3\CPb$, $3\CP\#5\CPb$ and $3\CP\#7\CPb$. The very same idea will apply to the others, as we
will discuss briefly. We begin by describing these families of 4-manifolds, showing that they all have the same
homeomorphism type, and afterwards we will use the Seiberg-Witten invariants to distinguish their diffeomorphism types. The SW invariants can be distinguished either via Theorem~\ref{thm: FS infinite} or using techniques in \cite{FS: same SW}.

In either method of computation, we first need to choose a null-homologous torus and peform $1/n$\/ surgery on it as in
Subsection~\ref{subsec: FS infinite}. We then prove that $\pi_1 = 1$ for the resulting infinite family of $4$-manifolds. To apply Theorem~\ref{thm: FS infinite} in its full strength, i.e. to obtain a family that consists of pairwise nondiffeomorphic 4-manifolds, we will show that we have exactly one basic class for $X_{\Lambda,\lambda}(0)$, up to sign, for each exotic $X$\/ that we have constructed. The authors believe that this should essentially follow from the fact that in all our constructions we obtain the manifold $X'$ after fiber summing two symplectic manifolds along a genus two surface. A folklore conjecture is that any manifold obtained this way should have the canonical class, up to sign, as its only Seiberg-Witten basic class (see for example \cite{MW}).  We will do this check by straightforward calculations using adjunction inequalities.

In all the constructions in Section~\ref{sec:general}, we observe that there is a copy of $V=(S^1 \times M_K) \setminus (F \cup S)$ embedded in the exotic $X$\/ we constructed, where $F$\/ is the fiber and $S$\/ is the section of $S^1\times M_K$, viewed as a torus bundle over a torus.  As shown in the Appendix, $S^1 \times M_K$ is obtained from $T^4$ after two Luttinger surgeries, which are performed in the complement of $F \cup S$. So we can think of $S^1 \times M_K$ as being obtained in two steps. Let $V_0$ be the complement of $F \cup S$ in the intermediate 4-manifold which is obtained from $T^4$ after the first Luttinger surgery. The next Luttinger surgery, say $(L,\gamma,-1)$, produces $V$\/ from $V_0$.  (In the Appendix, $L=c\times\tilde{b}$ and $\gamma=\tilde{b}$.)  This second  surgery on $L$\/ in $V_0$ gives rise to a nullhomologous torus $\Lambda$ in $V$. There is a loop $\lambda$ on $\Lambda$ so that the $0$ surgery on $\Lambda$ with respect to $\lambda$\/ gives $V_0$ back.  As the framing for this surgery must be the nullhomologous framing, we call it the `$0$-framing'.  Note that performing a $1/n$\/ surgery on $\Lambda$ with respect to $\lambda$ and this 0-framing in $V$\/ is the same as performing an $(L, \gamma, -(n+1))$ surgery in $V_0$ with respect to the Lagrangian framing.  We denote the result of such a surgery by $V(n)=V_{\Lambda,\lambda}(1/n)$. In this notation, $V(\infty)=V_0$ and we see that $V(0) = V$.

From the Appendix, we know that performing a $-n$ surgery on $L$\/ with respect to $\gamma$\/ and the Lagrangian framing, we obtain $V(n-1) = (S^1 \x M_{K_n}) \setminus (F \cup S)$, where $K_n$ is the $n$-twist knot.
It should now be clear that replacing a copy of $V$\/ in $X$\/ with $V(n-1) = (S^1 \x M_{K_n}) \setminus (F \cup S)$ (i.e.\ `using the $n$-twist knot') has the same effect as performing a $1/(n-1)$\/ surgery in the 0-framing on $\Lambda$ in $V \subset X$. We denote the result of such a surgery by $X_n=X_{\Lambda,\lambda}(1/(n-1))$. Clearly, $X_1=X$.  We claim that the family $\{ X_n \mid n=1,2,3,\dots \}$ are all homeomorphic to $X$\/ but have pairwise inequivalent Seiberg-Witten invariants. The first claim is proved in the following lemma.

\begin{lem}\label{lemma:pi_1(X_n)}
Let $X_n$ be the infinite family corresponding to a fixed exotic copy of $\CP \# \, (2k+1) \CPb$ for $k=1, \ldots, 4$ or of \,$ 3 \/ \CP \# (2l+3) \CPb$ for $l =1, \ldots,  6$ that we have constructed above. Then $X_n$ are all homeomorphic to $X$.
\end{lem}

\begin{proof}
For a fixed exotic $X$, $\pi_1(X_n)$ only differs from $\pi_1(X)$ by replacing a single  relation of the form $b = [a^{-1},d]$ by $b = [a^{-1},d]^n$ in the presentation of $\pi_1(X)$ we have used. (See the Appendix.)  Thus one only needs to check that raising the power of the commutator in one such relation does not effect our calculation of $\pi_1(X)=1$. This is easily verified in all of our examples. Hence all the fundamental group calculations follow the same lines and result in the trivial group.

Since $X_n$'s differ from $X$\/ only by surgeries on a nullhomologous torus, the characteristic numbers remain the same. On the other hand, since none have new homology classes, the parity should be the same.  By Freedman's theorem again, they all should be homeomorphic to each other.
\end{proof}

\medskip
Below, let $X$\/ be a 4-manifold obtained by fiber summing $4$-manifolds $Y$\/ and $Y'$ along submanifolds
$\Sigma\subset Y$\/ and $\Sigma'\subset Y'$. Let $A \subset Y$\/ and $B' \subset Y'$ be surfaces transversely intersecting $\Sigma$ and $\Sigma'$ positively at one point, respectively. Then we can form the internal connected
sum $A \# B'$ inside the fiber sum $X$, which is the closed surface that is the union of punctured surfaces
$(A\setminus(A\cap\nu\Sigma))\subset (X \setminus \nu\Sigma)$ and $(B'\setminus(B'\cap\nu\Sigma'))\subset (Y' \setminus\nu\Sigma')$. It is not hard to see that the intersection number between $A\# B'$ and $\Sigma=\Sigma'$ in $X$\/ is one, and thus they are both homologically essential. If all these manifolds and submanifolds are symplectic and the fiber sum is done symplectically, then $A\# B'$ can be made a symplectic submanifold of $X$\/ as well. Also note that, if either $A$\/ or $B'$ has self-intersection zero, then their parallel copies in their tubular neighborhoods can also be used to produce such internal sums in $X$.

\medskip
\subsection{An infinite family of exotic $\CP\# 3\CPb$'s}
\label{subsec:infinite1-3}

Let $X$\/ be the exotic $\CP\# 3\CPb$ that we have described in Subsection~\ref{subsec:1-3}. We begin by describing the surfaces that generate $H_2(X;\Z)$. There is a symplectic torus $T=T^2\times\{{\rm pt}\}$ of self-intersection zero in $Y=T^2\times \Sigma_2$ intersecting $\Sigma_2=\{{\rm pt}\}\times\Sigma_2$ positively at one point. On the other side, in $Y'=T^4\#2\CPb$, there is a symplectic torus $T'_1$ of self-intersection zero, and two exceptional spheres $E'_1$ and $E'_2$, each of which intersects $\Sigma'$ positively at one point.
(There is actually another symplectic torus $T'_2$ in $Y'$ satisfying $[\Sigma']=[T'_1]+[T'_2]-[E'_1]-[E'_2]$ in $H_2(Y';\Z)$, but we will be able to express the homology class that $T'_2$ induces in $X$\/ in terms of the four homology classes below.)

Hence we have four homologically essential symplectic surfaces: two genus two surfaces $\Sigma_2 = \Sigma'$,
$G = T \# T'_1$, and two tori $R_i = T \# E_i'$, $i = 1, 2$. Clearly $[\Sigma_2]^2 = [G]^2 = 0$, and $[R_1]^2 =[R_2]^2= -1$. It is a straightforward argument to see that these span $H_2(X;\Z)$, and the corresponding intersection form is isomorphic to that of $\CP\# 3\CPb$.  (Note that $[T\#T'_2]=[\Sigma_2]-[G]+[R_1]+[R_2]$.)

The $0$-surgery on $\Lambda$ with respect to $\lambda$ results in a 4-manifold $X_0=X_{\Lambda,\lambda}(0)$ satisfying $H_1(X_0;\Z)\cong \Z$ and $H_2(X_0;\Z) \cong H_2(X;\Z) \oplus \Z^2$, where the new 2-dimensional homology classes are represented by two Lagrangian tori $L_1$ and $L_2$. Both $L_j$ have self-intersection zero. They intersect each other positively at one point, and they do not intersect with any other class. Thus the adjunction inequality forces this pair to not  appear in any basic class of $X_0$.  Denoting the homology classes in $X_0$ that come from $X$\/ by the same symbols, let $\beta = a [\Sigma_2] + b [G] + \sum_i r_i [R_i]$ be a basic class of $X_0$.  Since it is a characteristic
element, $a$\/ and $b$\/ should be even, and $r_1$ and $r_2$ should be odd.

Since $b_2^+(X_0) > 1$, applying the (generalized) adjunction inequality for Seiberg-Witten basic classes (cf.\ \cite{OS:adj ineq}) to all these surfaces, we conclude the following.
\begin{itemize}
\item[(i)] \ $2 \geq 0 + | \beta \cdot [G] |$, implying $2 \geq |a|$.
\item[(ii)] \ $2 \geq 0 + | \beta \cdot [\Sigma_2] |$, implying $2 \geq |b + \sum_i r_i |$.
\item[(iii)] \ $0 \geq -1 + | \beta \cdot [R_i] |$, implying
$ 1 \geq | a - r_i|$ for $i = 1, 2$.
\end{itemize}
On the other hand, since $X_0$ is symplectic and $b_2^+(X_0) > 1$, $X_0$ is of simple type so we have $\beta^2 = 2e(X_0) + 3 \sigma(X_0) = 2e(X) + 3\sigma(X) = 2\cdot 6 + 3(-2) =  6$, implying:
\begin{itemize}
\item[(iv)] \ $6 = 2 a (b + \sum_i r_i) - \sum_i r_i^2$.
\end{itemize}
From (i), we see that $a$\/ can only be $0$ or $\pm 2$. However, (iv) implies that $a \neq 0$. Let us take $a=2$. Then by (iv) and (ii) we have
\begin{equation*}
 \Big| 6 + \sum_{i} r_i^2 \Big| = 4 \Big| b + \sum_i r_i \Big| \leq 8 ,
\end{equation*}
which implies that $ \sum_i r_i^2 \leq 2$. Therefore by (iii) we see that both $r_i$ have to be $1$. Finally by (iv) again, $b = 0$.

Similarly, if we take $a= -2$, we must have $r_1 = r_2 = -1$ and $b=0$.  Hence the only basic classes of $X_0$ are $\pm(2 [\Sigma_2] + \sum_i [R_i])=\pm K_{X_0}$, where $K_{X_0}$ denotes the canonical class of $X_0$. By Theorem~\ref{thm: FS infinite}, all $X_n=X_{\Lambda,\lambda}(1/(n-1))$ are pairwise nondiffeomorphic.

Moreover, by Remark~\ref{rmk: two classes} we see that $X = X_1$ also has one basic class up to sign.  It is easy to see that this is the canonical class $K_X = 2 [\Sigma_2] + [R_1] + [R_2]$.  Therefore the square of the difference of the two basic classes is $4K_X^2=24\neq -4$, implying that $X$\/ is irreducible (and hence minimal) by a direct application of Seiberg-Witten theory (cf.\ \cite{FS: immersed spheres}).  Furthermore, the basic class $\beta_n$ of $X_n$ corresponding to the canonical class $K_{X_0}$ satisfies
\begin{eqnarray}\label{eq:SW of 1-3}
\SW_{X_n}(\beta_n) &=& \SW_{X}(K_X) + (n-1)\SW_{X_0}(K_{X_0})\\
&=&  1 + (n-1) = n . \nonumber
\end{eqnarray}
Thus every $X_n$ with $n\geq 2$\/ is nonsymplectic.  In conclusion, we have proved the following.

\begin{thm}\label{thm:infinite1-3}
There is an infinite family of pairwise nondiffeomorphic\/ $4$-manifolds which are all homeomorphic to $\CP \# 3 \CPb$.  All of these manifolds are irreducible, and they possess exactly one basic class, up to sign.  All except for one are nonsymplectic.
\end{thm}

\subsection{Infinite smooth structures on $\CP\#(2k+1)\CPb$ for $k = 2,3,4$}
\label{subsec:infinite1-(2k+1)}

For exotic $\CP\# 5\CPb$'s, the second homology of $X_0$ will be generated by the following surfaces: two genus two
surface of self-intersection zero, $\Sigma = \Sigma'$ and $G= T \# T'$, four tori of self-intersection $-1$, $R_i = E_i \# T'$ ($i=1,2$) and $S_j = T \# E'_j$ ($j=1,2$), and two Lagrangian tori $L_1$ and $L_2$ as before.  A basic class of $X_0$ is of the form $\beta = a [\Sigma] + b [G] +  \sum_i r_i [R_i] + \sum_j s_j [S_j]$, where $a$\/ and $b$\/ are even and $r_i$ and $s_j$ are odd. The inequalities are:
\begin{itemize}
\item[(i)]  \ $2 \geq 0 + | \beta \cdot [G] |$, implying $ 2 \geq |a|$.
\item[(ii)]  \ $2 \geq 0 + | \beta \cdot [\Sigma] |$,
implying  $2 \geq |b + \sum_i r_i  + \sum_j s_j |$.
\item[(iii)] \  $0 \geq -1 + | \beta \cdot [R_i] |$, implying
$1 \geq | a - r_i|$ for $i=1,2$.
\item[(iv)] \ $0 \geq -1 + | \beta \cdot [S_j] |$, implying
$1 \geq | a - s_j|$ for $j=1,2$.
\item[(v)]  \  $4 = 2 a (b + \sum_i r_i + \sum_j s_j) - ( \sum_i r_i^2 + \sum_j s_j^2)$.
\end{itemize}
By (i), $a$\/ can only take the values $0, \pm 2$, where $0$ is ruled out by looking at (v). If $a=2$, then by (iii) and (iv) $r_i$ and $s_j$ are either $1$ or $3$. However, using (ii) and (v) as before, we see that none of these can be $3$. It follows that $r_i = s_j =1$ for all $i$\/ and $j$, and $b= -2$. The case when $a = -2$ is similar, and we see that $X_0$ has only two basic classes $\pm (2 [\Sigma] - 2[G] + \sum_i [R_i] + \sum_j [S_j])$.

For exotic $\CP\# 7 \CPb$'s, the classes are similar except now we have four tori of the form $S_j$.  For a basic class $\beta = a [\Sigma] + b[G] +  \sum_{i=1}^2 r_i [R_i] + \sum_{j=1}^4 s_j [S_j]$ of $X_0$, we see that the coefficients have the same parity as above.  The first four inequalities are the same (with (iv) holding for $j=1,\dots,4$), whereas the last equality (v) coming from $X_0$ being of simple type becomes:
\begin{itemize}
\item[(v$'$)] $2 = 2 a (b + \sum_{i=1}^2 r_i + \sum_{j=1}^4 s_j)
- (\sum_{i=1}^2 r_i^2 + \sum_{j=1}^4 s_j^2)$.
\end{itemize}
Once again, (i) implies that $a$\/ is $0$ or $\pm 2$, but by (v$'$) it cannot be $0$. If $a=2$, then by exactly the same argument as before we see that $r_i = s_j =1$ for all $i$\/ and $j$, and thus $b= -4$. The case $a = -2$ is similar. Therefore the only two basic classes of $X_0$ are $\pm (2 [\Sigma] - 4 [G] + \sum_{i=1}^2 [R_i] + \sum_{j=1}^4 [S_j])$.

For exotic $\CP\# 9\CPb$'s, the only difference is that the number of tori in total is eight. Let us denote two of the additional tori as $R_3$ and $R_4$ corresponding to say $R_a$ and $\Sigma_a$ where the other two will be denoted by $S_3$ and $S_4$ corresponding to $R_b$ and $\Sigma_b$ as described in Subsection \ref{subsec:1-9}.

For a basic class $\beta = a [\Sigma] + b [G] +  \sum_{i=1}^4 r_i [R_i] + \sum_{j=1}^4 s_j [S_j]$, once again $a, b$ are even and $r_i, s_j$ are odd.  The inequalities (i)--(iv) remain the same. Finally (v$'$) becomes:
\begin{itemize}
\item[(v$''$)]  $0 = 2 a (b + \sum_{i=1}^4 r_i + \sum_{j=1}^4 s_j)
- ( \sum_{i=1}^4 r_i^2 + \sum_{j=1}^4 s_j^2)$.
\end{itemize}
As before, $a$\/ cannot be $0$.  If $a =2$, then by the same argument we see that $r_i = s_j =1$ for all $i$\/ and $j$.  Thus $b= -6$, and we get a basic class $\beta = 2 [\Sigma] -6 [G] + \sum_{i=1}^4 [R_i] + \sum_{j=1}^4 [S_j]$. For $a =-2$ it is easy to check that we get the negative of this class.

Hence in all three examples $X_0$ has only two basic classes, and therefore by Theorem~\ref{thm: FS infinite},
all three families $\{X_n\}$ consist of pairwise nondiffeomorphic 4-manifolds homeomorphic to $\CP \# 5 \CPb$, $\CP \# 7 \CPb$, or $\CP \# 9 \CPb$. Furthermore, as in the previous subsection, we see that each family $\{X_n\}$ consists of 4-manifolds with only one basic class, up to sign. In each family, all but one member are nonsymplectic as the only nonzero values of the Seiberg-Witten invariant of $X_n$ are $\pm n$.  Finally, each exotic $X_n$\/ can be seen to be irreducible by a direct Seiberg-Witten argument as before (cf.\ \cite{FS: immersed spheres}).

\medskip
\subsection{An infinite family of exotic $3\CP\# 5\CPb$'s}
\label{subsec:infinite3-5}

Now let $X$\/ denote the exotic $3\CP\# 5\CPb$ that we have described in Subsection~\ref{subsec:3-5}.  Recall that we blow up $T^4$ twice and fiber sum along genus two surfaces of self-intersection $0$ in $Y=\Sigma_2\times\Sigma_2$ and in $Y' = T^4\#2\CPb$.

We begin by describing the surfaces that generate $H_2(X;\Z)$. There is a symplectic genus two surface $D$\/ of self-intersection zero in $Y$\/ intersecting $\Sigma$ positively at one point. On the other side in $T^4 \# 2 \CPb$ there is a symplectic torus $T'$ of self-intersection zero, intersecting $\Sigma'$ positively at one point. (There is also another symplectic torus orthogonal to this one, but the homology class it induces in $X$\/ can be expressed in terms of the other homology classes that we will present below.) There are also two exceptional spheres $E_1'$ and $E_2'$ in $T^4 \# 2 \CPb$, each of which intersects $\Sigma'$ positively at one point as well. In addition to these there are four Lagrangian tori $L_j\subset Y'$, $j=1, \dots, 4$, which come in hyperbolic pairs, and which do not intersect $\Sigma'$. Thus these four tori are certainly contained in $X$. The other four homology generators in $X$\/ are represented by the following symplectic surfaces: a genus two surface $\Sigma = \Sigma'$, a genus three surface $G = D\# T'$, and genus two surfaces $R_i = D \# E_i'$, $i = 1, 2$. We have $[\Sigma]^2 = [G]^2 = 0$ and $[R_i]^2 = -1$. Once again it is a straightforward calculation to show that these eight elements induce an intersection form which is isomorphic to that of $3\CP\# 5\CPb$.

The $0$ surgery on $\Lambda$ results in a 4-manifold $X_0$ satisfying $H_1(X_0;\Z)= \Z$ and $H_2(X_0;\Z) \cong H_2(X) \oplus \Z^2$, where the new 2-dimensional homology classes are represented by two Lagrangian tori $L_1^0$ and $L_2^0$.  As we have argued in Subsection~\ref{subsec:infinite1-3}, the homology classes $[L_j^0]$ do not appear in any basic
class of $X_0$, and neither do the other four Lagrangian tori $L_j$. Thus a basic class of $X_0$ is of the form $\beta = a [\Sigma] + b [G] + \sum_i r_i [R_i]$ as before. Since it is a characteristic element, $a$ and $b$\/ should be
even, and $r_1$ and $r_2$ should be odd.

The inequalities (i)--(iv) in Subsection~\ref{subsec:infinite1-3} now become (due to the change in genera of the surfaces and the characteristic numbers of $X_0$):
\begin{itemize}
\item[(i)] \ $4 \geq 0 + | \beta \cdot [G] |$, implying $4 \geq |a|$.
\item[(ii)] \ $2 \geq 0 + | \beta \cdot [\Sigma] |$, implying $2 \geq |b + \sum_i r_i |$.
\item[(iii)] \ $2 \geq -1 + | \beta \cdot [R_i] |$, implying $3 \geq | a - r_i|$
for $i = 1, 2$.
\item[(iv)] \ $14 = 2 a (b + \sum_i r_i) - \sum_i r_i^2$.
\end{itemize}
From (i), we see that $a$ can only be $0,\pm 2, \pm 4$. However, (iv) indicates that $a \neq 0$. Let us take $a=2$. Then by (iv) $(b + \sum_i r_i) \geq 4$, which contradicts (ii). Similarly if $a= -2$, we get $(b + \sum_i r_i) \leq - 4$, again contradicting (ii).

Next suppose that $a = 4$. Then by (iv) $(b + \sum_i r_i) \geq 2$ and so (ii) implies that $(b+ \sum_i r_i) =2$.  Now (iv) implies that $\sum_i r_i^2=2$, and hence $r_1 = r_2 = \pm 1$. By (iii) we see that $r_1=r_2= 1$.  Then $b = 0$.  Therefore we get $\beta = 4[\Sigma] + \sum_i [R_i]$. The case when $a = -4$ is similar. As before we see that $r_1 = r_2 = \pm 1$ and (iii) forces both $r_i$ to be equal to $-1$. It again follows that $b = 0$. Hence $\beta$\/ is the negative of the class we obtained for the $a=4$ case.

It follows that $X_0$ has only two basic classes $\pm (4 [\Sigma] + \sum_i [R_i])=\pm  K_{X_0}$, where $K_{X_0}$ is the canonical class of the symplectic manifold $X_0$. Once again by Theorem~\ref{thm: FS infinite}, all $X_n$ are pairwise nondiffeomorphic. So we obtain the following.

\begin{thm}
There is an infinite family of pairwise nondiffeomorphic\/ $4$-manifolds which are all homeomorphic to $3\CP \# 5 \CPb$.  All of these\/ $4$-manifolds are irreducible, and they possess exactly one basic class, up to sign.  All except for one are nonsymplectic.
\end{thm}

\medskip
\subsection{Infinite smooth structures on $3\CP\# (2l+3)\CPb$ for $l = 2,\dots, 6$}
\label{subsec:infinite3-(2l+3)}

In Subsection~\ref{subsec:3-(2l+3)} we constructed irreducible symplectic exotic $3\CP\# (2l+3)\CPb$ for $l = 2,\dots, 6$ by fiber summing irreducible symplectic exotic $\CP \# (2k+1) \CPb$ for $k =1, \ldots, 4$ along a genus two surface with $T^4 \# 2 \CPb$ or $(S^2 \x T^2) \# 4 \CPb$.  For each $3\CP\# (2l+3)\CPb$ for $l = 2,\dots, 6$, we will obtain an infinite family of exotic smooth structures on it, by simply using the infinite family of exotic smooth structures on $\CP \# (2k+1) \CPb$ constructed in Subsections \ref{subsec:infinite1-3} and \ref{subsec:infinite1-(2k+1)}.  Since each such infinite family was obtained by performing torus surgeries away from the genus two surface $\Sigma$ that we have used to fiber sum, we can fiber sum each exotic member along $\Sigma$ with the usual genus two surface in $T^4 \# 2 \CPb$ or $(S^2 \x T^2) \# 4 \CPb$.  These new fiber sums are then also simply-connected, and they have the same characteristic numbers and parity as the exotic copy of $3\CP\# (2l+3)\CPb$ ($l = 2,\dots, 6$) that we have constructed. We claim that infinitely many of them have different Seiberg-Witten invariants. We will present the details of this argument for the family which is homeomorphic to $3 \CP \# 7 \CPb$ and leave the others as an exercise for the reader since they are extremely similar.

\begin{thm}\label{thm:infinite3-7}
There is an infinite family of pairwise nondiffeomorphic\/ $4$-manifolds which are all homeomorphic to\/ $3\CP \# 7 \CPb$.
These\/ $4$-manifolds are all irreducible, and all except for one are nonsymplectic.
\end{thm}

\begin{proof}
Let $\{X_n\mid n = 1,2,3,\dots \}$ be the infinite family of exotic $\CP\#3\CPb$ in Lemma~\ref{lemma:pi_1(X_n)} and Theorem~\ref{thm:infinite1-3}.  We have shown in Subsection~\ref{subsec:infinite1-3} that each $X_n$ has only one basic class up to sign (cf.\ Remark~\ref{rmk: two classes}).  Let us denote these basic classes by $\pm\beta_n\in H_2(X_n;\Z)$.  From (\ref{eq:SW of 1-3}) in Subsection~\ref{subsec:infinite1-3}, we conclude that
\begin{equation}\label{eq:SW of X_n}
\SW_{X_n}(\pm\beta_n)=\pm n.
\end{equation}
Using the adjunction inequality and the blowup formula in \cite{FS: immersed spheres}, we can deduce that there are four basic classes of $T^4\#2\CPb$, namely
$[E_1]\pm [E_2]$ and $-[E_1]\pm[E_2]$, where $E_1$ and $E_2$ are the exceptional spheres of the blowups.  We can choose orientations so that the values of $\SW_{T^4\#2\CPb}$ on all four basic classes are $1$.  We need to compare the SW invariants of the fiber sums $Z'_n=X_n\#_{\Sigma_2=\Sigma_2'}(T^4\#2\CPb)$.  These are non-symplectic fiber sums when $n>1$.   Note that we can view $Z'_n$ as the result of $1/(n-1)$ surgery on a nullhomologous torus in the symplectic sum $Z'=X\#_{\Sigma_2=\Sigma_2'}(T^4\#2\CPb)$ that was constructed in Corollary~\ref{cor: 3CPs}.  There is a canonical isomorphism $\alpha_n:  H_2(Z';\Z)\to H_2(Z'_n;\Z)$.  Using (\ref{eq:SW of X_n}) and the product formula in \cite{MST}, we can easily compute that
\begin{equation*}
\SW_{Z'_n}(\pm \kappa) = \pm n,
\end{equation*}
where $\kappa=\alpha_n({\rm PD}(K_{Z'}))$ and ${\rm PD}(K_{Z'})$ is the Poincar\'{e} dual of the canonical class of $Z'$.  Since $Z'$ is irreducible by Corollary~\ref{cor: 3CPs},  $Z'_n$ has to be irreducible as well.
\end{proof}

\section{Open questions}
\label{sec:questions}

We would like to finish with two questions that naturally arise within the context of our article.
\subsection*{Question 1} Is there a general construction scheme similar to the one we have outlined, which produces irreducible symplectic structures on $\CP \# 2k \CPb$ or on $3 \CP \# 2l \CPb$ (for $0 \leq k \leq 4$ and for $0 \leq l \leq 5$)? More explicitly, what would be the intermediate manifold $X'$, satisfying conditions (I) and (II)
in Section~\ref{sec:general}, if the target manifold $Z$\/ is chosen as one of these manifolds?

Certainly, for $k \leq 2$ and $l \leq 3$, it is currently unknown whether there are any irreducible $4$-manifolds in these homeomorphism classes.  When $k=3$ or $4$, \cite{SS} and \cite{Pa:2005} give irreducible symplectic structures on $\CP\#2k\CPb$ and \cite{FS: double node} give infinite families of pairwise nondiffeomorphic exotic smooth structures.
For $l=4$ or $5$, we refer to \cite{Pa:3-8} and \cite{P3}. Our emphasis here is on the method, and it would be nice to construct similar examples using our recipe in Section~\ref{sec:general}.

On the other hand, note that Remark \ref{the rest} can easily be modified to obtain exotic $3 \CP \# (2k + 10) \CPb$ for $k = 1, \ldots, 4$. For this, we first blow down one of the sphere sections in $E(1)$ to get a symplectic torus $T_2$ of self-intersection $+1$ in $\CP \# 8 \CPb$, and then fiber sum $\CP \# 8 \CPb$ with an exotic copy $X$\/ of $\CP \# (2k+1) \CPb$ for $k = 1, \ldots, 4$ by gluing $T_2$ to one of the symplectic tori of self-intersection $-1$ in $X$.

\medskip
Our next question is motivated by our $3 \CP \# 5 \CP$ example as well as the constructions that appear in \cite{FPS}.

\subsection*{Question 2} Is it possible to determine the conditions on the fundamental group of the intermediate manifold $X'$, possibly with some geometric side-conditions, that allow one to kill the fundamental group through a sequence of Luttinger surgeries?

\medskip
\section{Appendix}
\label{sec:appendix}

\begin{figure}[ht]
\begin{center}
\includegraphics[scale=.7]{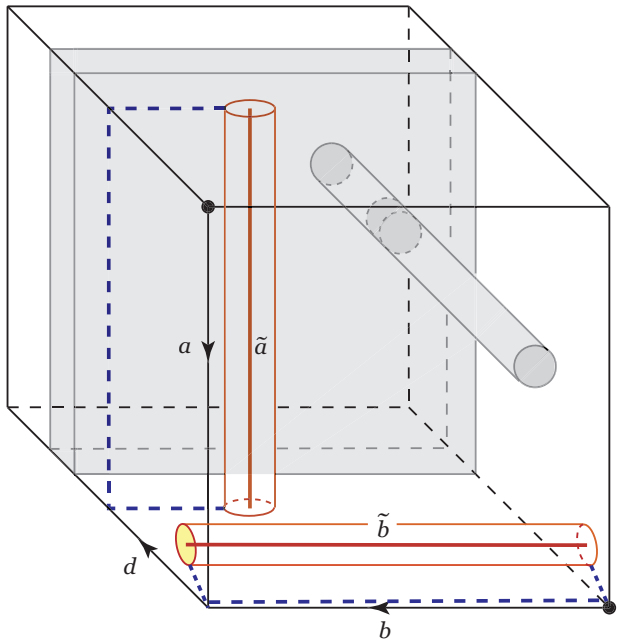}
\caption{The 3-torus $d\times a \times b$}
\label{fig:cube}
\end{center}
\end{figure}

\begin{figure}[ht]
\begin{center}
\includegraphics[scale=.9]{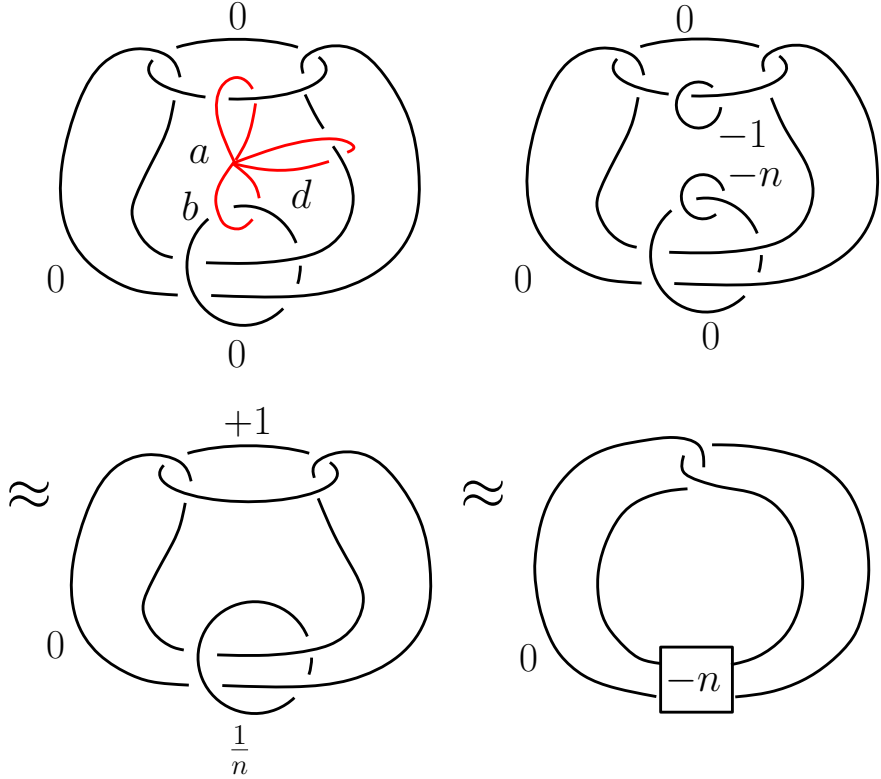}
\caption{\small The first diagram depicts the three loops $a, b, d$\/ that generate the
$\pi_1(T^3)$. The curves $\tilde{a} = d a d^{-1}$ and $\tilde{b}$\/ are freely
homotopic to the two extra curves given in the second diagram. The third
diagram is obtained from the second via two slam-dunk operations; wheras the
last diagram is obtained after Rolfsen twists.}
\label{fourtorus2}
\end{center}
\end{figure}

Let $T^4 =a \times b \times c \times d \cong (c \times d) \times (a \times b)$, where we have switched the order of the symplectic $T^2$ components $a\times b$\/ and $c\times d$\/ to have a comparable notation for our earlier $\pi_1$ calculations.  Let $K_n$ be an $n$-twist knot (cf. Figure~\ref{fourtorus2}). Let $M_{K_n}$ denote the result of performing $0$ Dehn surgery on $S^3$ along $K_n$.  Our goal here is to show that the 4-manifold $S^1 \times M_{K_n}$ is obtained from $T^4 = (c \times d) \times (a \times b) = c \times ( d \times a \times b) = S^1 \times T^3$ by first performing a Luttinger surgery $(c \times \tilde{a}, \tilde{a}, -1)$ followed by a surgery $(c \times \tilde{b}, \tilde{b}, -n)$.   Here, the tori $c\times \tilde{a}$\/ and $c\times\tilde{b}$\/ are Lagrangian and the second tilde circle factors in $T^3$ are as pictured in Figure~\ref{fig:cube}.   We use the Lagrangian framing to trivialize their tubular neighborhoods, so when $n=1$ the second surgery is also a Luttinger surgery.

Note that the normal disks of each Lagrangian tori in their Weinstein neighborhoods lie completely in $T^3$ and are disjoint. Thus topologically, the result of these surgeries can be seen as the product of the first $S^1$ factor with the result of Dehn surgeries along $\tilde{a}$ and $\tilde{b}$ in $T^3$.  Therefore we can restrict our attention to the effect of these Dehn surgeries in $T^3$ since the diffeomorphisms of the $3$-manifolds induce diffeomorphisms between the product $4$-manifolds.

The Kirby calculus diagrams in Figure~\ref{fourtorus2} show that the result of these Dehn surgeries is the manifold $M_{K_n}$, where $K_n$ is (the mirror of) the $n$-twist knot.  In particular, note that for $n=1$ we get the trefoil knot $K$. Thus the effect of $(c \times \tilde{b}, \tilde{b}, -n)$ surgery with $n>1$ as opposed to the Luttinger surgery $(c \times \tilde{b}, \tilde{b}, -1)$ is equivalent to using the non-symplectic 4-manifold $S^1 \times M_{K_n}$ instead of symplectic $S^1 \x M_K$ in our symplectic sum constructions.

Next we describe the effect of these surgeries on $\pi_1$.
First it is useful to view $T^3=d\times(a\times b)$\/ as a $T^2$ bundle over $S^1$
with fibers given by $\{{\rm pt}\}\times(a\times b)$ and sections given by $d\times\{{\rm pt}\}$.  The complement of a fiber union a section in $T^3$ is the complement of
3-dimensional shaded regions in Figure~\ref{fig:cube}.

\begin{figure}[ht]
\begin{center}
\includegraphics[scale=.5]{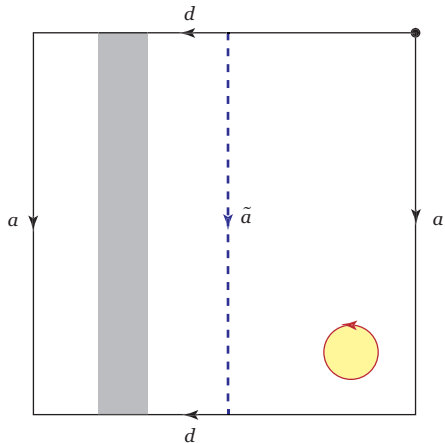}
\caption{\small The face of cube where we can see the Lagrangian pushoff of $\tilde{a}$}
\label{fig:cubeface}
\end{center}
\end{figure}

It is not too hard to see that the Lagrangian framings give the following product decompositions of two boundary 3-tori (compare with \cite{BK:1-3,FPS}):
\begin{eqnarray}
\partial(\nu(c\times\tilde{a})) &\cong& c \times (dad^{-1}) \times [d,b^{-1}],
\label{eq:c x a}\\
\partial(\nu(c\times\tilde{b})) &\cong& c \times b \times [a^{-1},d]. \label{eq:c x b}
\end{eqnarray}

The Lagrangian pushoff of $\tilde{b}$\/ is represented by $b$, as a homotopy to $b$\/ is given by the ``diagonal'' path (dotted lines emanating from the horizontal boundary cylinder $\partial(\nu\tilde{b})$ in Figure~\ref{fig:cube}).  For decomposition (\ref{eq:c x b}), it is helpful to view the base point as the front lower right corner of the cube represented by a dot in Figure~\ref{fig:cube}.  It is comparatively more difficult to see that the Lagrangian pushoff of $\tilde{a}$ is represented by $dad^{-1}$. The Lagrangian pushoff of $\tilde{a}$\/ is represented by the dotted circle in Figure~\ref{fig:cubeface} and is seen to be homotopic to the composition $a[a^{-1},d]=a(a^{-1}dad^{-1})=dad^{-1}$.  For decomposition (\ref{eq:c x a}), it is helpful to view the base point as the front upper left corner of the cube represented by a dot in Figure~\ref{fig:cube}.  The new relations in $\pi_1$ introduced by the two surgeries are
\begin{equation}\label{pi_1 relations from two surgeries}
dad^{-1} = [d,b^{-1}] =db^{-1}d^{-1}b, \quad b = [a^{-1},d]^n=(a^{-1}dad^{-1})^n .
\end{equation}

From now on, let us assume that $n=1$.
Then the second relation in (\ref{pi_1 relations from two surgeries}) gives
\begin{equation}\label{monodromy1}
ab = dad^{-1}.
\end{equation}
Combining (\ref{monodromy1}) with the first relation in (\ref{pi_1 relations from two surgeries}) gives $ab= dad^{-1} = db^{-1}d^{-1}b$, which can be simplified to
$a=db^{-1}d^{-1}$.  Thus we have
\begin{equation}\label{monodromy2}
a^{-1}=dbd^{-1}.
\end{equation}
Hence we see that (\ref{monodromy1}) and (\ref{monodromy2}) give the standard representation of the monodromy of the $T^2=a\times b$\/ bundle over $S^1=d$\/ that is the 0-surgery on $S^3$ along the trefoil $K=K_1$.

\medskip
\section*{Acknowledgments}

The authors thank Ronald Fintushel and Ronald J. Stern for helpful comments and discussions. A. Akhmedov was partially
supported by the NSF grant FRG-0244663, R. \.{I}. Baykur was partially supported by the NSF Grant DMS0305818, and B. D. Park was partially supported by CFI, NSERC and OIT grants.

\end{document}